\documentclass[conference]{IEEEtran}
\IEEEoverridecommandlockouts
\usepackage{cite}
\usepackage{amsmath,amssymb,amsfonts}
\usepackage{algorithmic}
\usepackage{graphicx}
\usepackage{textcomp}
\usepackage{xcolor}

\usepackage{comment}
\usepackage{color}
\usepackage{mathtools}
\usepackage{siunitx}
\usepackage{hyperref}

\usepackage{subcaption}
\usepackage{microtype}

\def\pd{\partial}

\begin{document}

\onecolumn
\begin{center}
   \Large
This work has been submitted to the I3DA 2021 International Conference for possible publication in the IEEE Xplore Digital Library. 
\\
\vspace{2em}
Copyright may be transferred without notice, after which this version may no longer be accessible.
\end{center}
\clearpage
\twocolumn
\title{Adding air attenuation to simulated room impulse responses: A modal approach}

\author{\IEEEauthorblockN{Brian Hamilton}
\IEEEauthorblockA{\textit{Acoustics \& Audio Group} \\
\textit{University of Edinburgh}\\
Edinburgh, UK \\
bhamilton.ac@gmail.com}
}

\maketitle

\begin{abstract}
   Air absorption is an important effect to consider when simulating room acoustics as it leads
   to significant attenuation in high frequencies.  In this study, an offline method
   for adding air absorption to simulated room impulse responses is devised.  The proposed
   method is based on a modal scheme for a system of one-dimensional dissipative wave equations,
   which can be used to post-process a room impulse response simulated without air absorption,
   thereby incorporating missing frequency-dependent distance-based air attenuation.  Numerical
   examples are presented to evaluate the proposed method, along with comparisons to existing
   filter-based methods.
\end{abstract}


\section{Introduction}
Air absorption is an important source of attenuation in room acoustics and its emulation is
critical to achieving realistic room acoustic simulations~\cite{vorlander2020auralization}.  In
general, air attenuation is a frequency-dependent and distance-based effect caused by
viscothermal effects and relaxation processes~\cite{pierce2019acousticsChap9}.  While these
physical processes may be modelled directly with wave-based
models~\cite{wochner2005atmospheric,jimenez2016time,hamilton2020viscothermal}, the effect of air
absorption is typically  modelled with a set of digital
filters~\cite{moorer1979reverberation,jot1991digital,huopaniemi1997modeling,savioja1999creating,schroder2011physically,saarelma2018challenges,kates2020airabsorption}.
Filter approaches have the potential to accurately reproduce air attenuation, but a common
underlying problem is that filters must be tuned to accurately apply air attenuation.  This
comes with a trade-off between simplicity and accuracy -- ranging from, e.g., 
octave-band approaches~\cite{scheibler2018pyroomacoustics} to optimised IIR
filters~\cite{kates2020airabsorption} and window-method FIR
filters~\cite{saarelma2018challenges}.  

Recently, a more physically-motivated filtering approach was proposed~\cite{hamilton2021dafx},
which is based on an approximate Green's function to the viscothermal wave equation and
side-steps many of the issues of pre-existing filter approaches.  However, that approximate
Green's function approach is best suited to controlled indoor air conditions where power-law
attenuation is
sufficient for audible frequencies.  The aim of this paper is to present a more general method
that allows for post-processing a pre-computed room impulse response (RIR) with a complete air
attenuation profile (i.e., classical power-law absorption plus two relaxation
effects~\cite{pierce2019acousticsChap9,iso19939613,bass1995atmospheric}) that can be
used for all air conditions (temperature and humidity).  This is particularly important when
simulating partly-outdoor spaces, such as, e.g., medieval cathedrals~\cite{postma2016acoustics}
or other heritage/historic
sites~\cite{iannace2014acoustics,murphy2017acoustic,katz2020exploring,dorazio2020understanding}.

The general aim of the proposed method is to use an input signal to drive a set of
exponentially-damped plane waves, which travel a given distance such that the resulting output
can be seen as frequency- and distance-based air attenuation having been applied to the input
signal.  As will be seen, this can be achieved with a system of dissipative wave equations which
can be solved accurately and efficiently by a modal time-stepping scheme.

This paper is structured as follows.  Section~\ref{secback} presents background theory for air
absorption and lossy wave propagation. Section~\ref{secnum1} gives numerical schemes for lossy
wave propagation, building towards a modal scheme of interest, along with numerical evaluations.   Section~\ref{secmodal} presents
the frequency-dependent modal scheme and procedure for adding air absorption to RIRs, along with
numerical experiments testing its performance against a reference solution.  
Conclusions and final remarks are given in Section~\ref{secconc}.

\begin{figure}[t]
   \centering
   \includegraphics[scale=0.48,clip,trim=1.8cm 7cm 2cm 7cm]{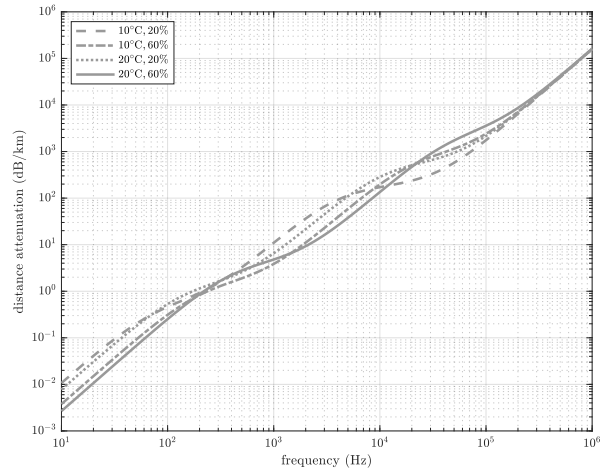}
   \caption{Distance attenuation in dB/km as a function of frequency, for varying air conditions
   (temperature in degrees Celsius and relative humidity in percent).}
   \label{figabs}
\end{figure}

\section{Background\label{secback}}
\subsection{Air absorption}
Air absorption can be seen as a frequency-dependent distance-based damping on plane waves, written as:
\begin{equation}
   p(t,x) =  \hat{p}e^{-\alpha(\omega)x}e^{i(k x -\omega t)} 
   \label{eqnplanewavespace}
\end{equation}
where $p$ is a pressure, with a complex amplitude $\hat{p}$, real wavenumber $k$, real angular frequency $\omega$, 
time $t$, and position~$x$.  $\alpha(\omega)$ is a
frequency-dependent attenuation coefficient in Nepers/m which represents air absorption.
Formulae for $\alpha(\omega)$ as a function of air temperature and humidity are left out for brevity but may be found
in~\cite{iso19939613,bass1995atmospheric}. See Fig.~\ref{figabs} for typical behaviour.
For the purpose of this study it is sufficient to characterise this absorption by the property that $\alpha(\omega)\ll \omega/c$,
such that $k\approx\omega/c$, where $c$ is the speed of sound in air, which also means that the
phase velocity can be assumed constant and equal to $c$.


\subsection{Dissipative wave equation}
In order to build up a physical system to achieve the desired result, we start by
considering the one-dimensional lossy wave equation:
\begin{equation}
   \pd_t^2 p  + 2\sigma \pd_t p - c^2 \pd_x^2 p = 0
   \label{eqnlossywave}
\end{equation}
where $p=p(t,x)$ is our wave-variable of interest (e.g., a pressure quantity), $\pd_t^2$ and $\pd_x^2$ denote
second-order partial derivatives with respect to time and space, respectively, and $\sigma$ is a loss constant in s$^{-1}$.
A dispersion relation can be derived by considering a trial solution of the
form: 
\begin{equation}
   p(t,x) =  \hat{p}e^{i(\hat{k} x -\omega t)}  = \hat{p}e^{-\alpha x}e^{i(k x -\omega t)} 
   \label{eqnplanealphaconst}
\end{equation}
with complex wavenumber $\hat{k}=k+i \alpha$, where $\alpha$ is a constant (i.e., so far
dissipation is not frequency-dependent).
By insertion into~\eqref{eqnlossywave} we get:
\begin{equation}
   \hat{k}^2 = \frac{\omega^2}{c^2} + \frac{2i\sigma\omega}{c^2} \quad \Rightarrow \quad \hat{k} = \frac{\omega}{c}\sqrt{\left(1+ \frac{2i\sigma}{\omega}\right)}
   \label{eqnlossydisp}
\end{equation}
We assume $\sigma\ll \omega$, for which we get $k\approx \omega/c$ and the approximate square root:
\begin{equation}
   \hat{k} \approx \frac{\omega}{c} + i\frac{\sigma}{c}
   \label{eqnlossydispapprox}
\end{equation}
Equating this to $\hat{k}=k+i \alpha$ we recover $\alpha\approx \sigma/c$.  

It will be useful to set a boundary condition for this dissipative wave system that induces
one-way travelling \mbox{(plane-)}wave solutions.  One manner this can be achieved is through a forced
boundary condition~\cite{kinsler2000fundamentals}.\footnote{One-way wave
equations exist and could also be used (e.g., convection-diffusion
equations~\cite{strikwerda2004finite}).  Second-order (two-way) wave
equations are preferred here for favourable properties of associated numerical schemes.}

Returning to the time-space
domain, let us consider this system on the domain $x\in[0,+\infty)$.  At the left
boundary we impose the condition:
\begin{equation}
   p(t,0) = u(t)
   \label{eqnbc0}
\end{equation}
where $u(t)$ is some input signal used to drive the signal (with, otherwise, zero initial
conditions for $x>0$).  

In the lossless case, $\sigma=0$, it is straightforward to check that this boundary
condition results in a rightward-travelling wave solution:
\begin{equation}
   p(t,x) =  u(t-x/c)
   \label{eqntrav}
\end{equation}

When $\sigma>0$, the system is dispersive and travelling-wave solutions -- analogous
to~\eqref{eqntrav} -- are not available.  However, in the case of a time-harmonic
$u(t)=\hat{u}e^{-j\omega t}$, it is straightforward to see that the solution on the domain
$x\in[0,\infty)$ takes the
form:
\begin{equation}
   p(t,x) =  \hat{u}e^{i(\hat{k} x -\omega t)} 
\end{equation}
with $\hat{k}$ given by~\eqref{eqnlossydisp}.  

Let us then assume $u(t)$ can be written
as the following sum of time-harmonic signals~$u_q(t)$:
\begin{equation}
   u(t) = \sum_{q=0}^{Q-1} u_q(t)\,, \quad u_q(t) = \hat{u}_q e^{-i\omega_q t} 
   \label{eqnufourier}
\end{equation}
where $\omega_q$ are distinct (real) angular frequencies.
Clearly, each component $u_q(t)$ would have the associated pressure solution $p_q(t,x) =
\hat{p}_q e^{-\alpha x}e^{j(k_qx-\omega_q t)}$, and by  
linearity we have the solution:
\begin{equation}
   p(t,x) = \sum_{q=0}^{Q-1} p_q(t,x) = \sum_{q=0}^{Q-1} \hat{u}_q e^{-\alpha x}e^{i(k_qx-\omega_q t)}
   \label{eqnsumpressures}
\end{equation}
Thus, the general solution to this system, for $x\in[0,\infty)$, may be seen as a Fourier decomposition of damped,
travelling plane waves, driven by the input signal $u(t)$.

\section{Numerical schemes for lossy wave equation\label{secnum1}}
In this section we consider numerical schemes to simulate this lossy wave equation system.
Consider a grid function $p^n_l \approxeq p(nT,(l+1/2)X)$ with grid spacing $X$ and time-step
$T$, chosen to be
$T=X/c$.  Furthermore let $u^{n}=u((n-1/2) T)$.  In terms of these discrete functions, the lossless
travelling wave solution~\eqref{eqntrav} would translate to:
\begin{equation}
   \label{eqnwavedisc}
   p^{n}_l = u^{n-l}
\end{equation}

\subsection{FDTD schemes}
Although the intention is to derive and use a modal scheme, it will be informative to first discuss a
simple FDTD scheme.  The simplest finite-difference time-domain (FDTD) scheme in this
setting~\cite{bilbao2009numerical} would be:
\begin{equation}
   \label{eqnfdtd}
   p^{n+1}_l = \frac{1}{1+\sigma T }\left(p^n_{l+1} +  p^n_{l-1} +(\sigma T-1) p^{n-1}_{l}\right) 
\end{equation}
This scheme is known to be stable for $\sigma\geq 0$~\cite{bilbao2009numerical}.  To implement the boundary condition
applied at $l=0$, we can simply force the value at $l=0$:
\begin{equation}
   \label{eqnbc1}
   p^{n+1}_0 =  u^{n+1}
\end{equation}
The above can be seen as a ``hard source'' -- as opposed to a ``soft source'', which would be
a forcing term that additive to a pressure
update~\cite{sheaffer2014sources,botts2014sources}.

In the lossless case ($\sigma=0$), this scheme permits the solution~\eqref{eqnwavedisc} and is thus
\emph{exact} (this can be verified by hand).  This implies the scheme~\eqref{eqnfdtd} is stable
in the lossless case.  It is expected that with the inclusion of the loss term $2\sigma \pd_t p$ --
discretised here using the trapezoidal rule -- stable is maintained under~\eqref{eqnbc1},
although a stability analysis with the ``hard source'' boundary is not immediately available.

A close alternative boundary condition featuring a ``soft source'' is also available as:
\begin{equation}
   \label{eqnbc2}
   p^{n+1}_0 =  \frac{1}{1+\sigma T}\left(p^n_{1} + p^n_{0} +(\sigma T-1) p^{n-1}_{0}\right) + F^n 
\end{equation}
where the ``soft source'' forcing term is denoted $F^n$:
\begin{equation}
   \label{eqnforce}
   F^n = \frac{1+\sigma T/2}{1+\sigma T} u^{n+1} - \frac{1-\sigma T/2}{1+\sigma T} u^n
\end{equation}
The derivation for this update appears in the Appendix.
It is straightforward to check that this scheme is also exact in the lossless case.
Furthermore, a discrete
energy balance maybe be obtained for this scheme (with $\sigma\geq 0$), from which stability can
be proven (this is left out for brevity).
In the presence of loss ($\sigma>0$) these schemes are no longer exact, but they remain
consistent with their underlying models.\footnote{This is primarily due to the
fact that the introduction of loss in~\eqref{eqnfdtd} upsets a delicate balance of approximation errors in space and
time from the lossless case~\cite{bilbao2009numerical}.}


\subsection{Semi-discrete modal schemes}
With the target of simulating the general case of frequency-dependent dissipation, it will be
convenient to consider modal schemes.  To start we consider a semi-discrete modal scheme for
this dissipative wave problem.  For this we take the view that the solution is composed of modes
(as per~\eqref{eqnsumpressures}) of the form:
\begin{equation}
p_q(t,x) = P_q(t)\cos(k_qx) 
\end{equation}
with $P_q(t)$ unknown.
Inserting this trial solution into our PDE gives:
\begin{equation}
\ddot{P}_q(t) + 2\sigma \dot{P}_q(t) - c^2k_q^2 P_q(t)  = 0
\end{equation}
This ODE system ($q=0,\dots,Q-1$) can be solved with leapfrog integration, on a temporal-grid function $P^n_q \approxeq
P_q(nT)$, or it can be solved with an 
``exact'' recursion for the damped harmonic oscillator~\cite{cieslinski2006simulations,bilbao2009numerical,botts2015extension}:
\begin{equation}
P_q^{n+1} = 2e^{-\sigma T}\cos(\omega_q T)P_q^{n} - e^{-2\sigma T}P_q^{n-1}
\label{eqnmodalupdate}
\end{equation}
where $\omega_q = \sqrt{c^2 k_q^2 + \sigma^2}$.  This recursion is ``exact'' for solutions of
the form $P^n_q = \hat{p}_q e^{-j\hat{\omega}_qnT}$, which can be verified by insertion into the
above.  In this semi-discrete time-space domain, the numerical solution would then be:
\begin{equation}
p(nT,x) = \sum_{q=0}^{Q-1} P^n_q \cos(k_q x)
\end{equation}

%

\subsection{Fully-discrete modal schemes}
At this point we consider a fully-discrete modal scheme for our wave problem.  For practical
purposes we must truncate the domain of interest to some length $L_x$, such that $x\in[0,L_x]$.
For our purposes, it is sufficient to assume $L_x = cT_d$ where $T_d$ is a simulation duration
of interest (this will be set to the duration of our pre-simulated input room impulse response).
We consider a grid function $p^n_l\approxeq p(nT,(l+1/2)X)$, with $l=0,\dots,N_x-1$, with $N_x =
\lceil L_x/X \rceil$, and $T=X/c$.  For the modes, we choose the orthonormal cosine basis:  
\begin{equation}
   \Phi_{q,l} = a_q\cos\left( \frac{q\pi(l+1/2)}{N_x}\right)
   \label{eqncosinePhi}
\end{equation}
where $a_0 = \sqrt{\frac{1}{N_x}}$ and $a_q = \sqrt{\frac{2}{N_x}}$ for $q>0$.
As such, the solution is now assumed to take the form:
\begin{equation}
p^n_l = \sum_{q=0}^{Q-1} P^n_q \Phi_{q,l}
\end{equation}
with $Q=N_x$.

It is insightful to write the full solution in matrix form.  Let us define:
\begin{subequations}
\begin{align}
\mathbf{p}^n &= \left[p^n_0,p^n_1,\dots,p^n_{N_x-1}\right]^{\mathrm{T}}
\\
\mathbf{P}^n &= \left[P^n_0,P^n_1,\dots,P^n_{Q-1}\right]^{\mathrm{T}}
\\
\mathbf{\Phi}_{q} &= \left[\Phi_{q,0},\Phi_{q,1},\dots,\Phi_{q,N_x-1}\right]^{\mathrm{T}}
\\
\mathbf{V} &= \left[\mathbf{\Phi}_{0} | \mathbf{\Phi}_{1}| \dots |\mathbf{\Phi}_{Q-1}\right]
\end{align}
\end{subequations}
where $\mathbf{V}$ is a $N_x \times N_x$ matrix, and $\mathbf{p}^n$, $\mathbf{P}^n$ and
$\mathbf{\Phi}_{q}$ are $N_x\times 1$ vectors (and ``$^\mathrm{T}$'' denotes transposition),
with $Q=N_x$.
With these we have the ability to transform between time-space and modal domains via the matrix-vector products:
\begin{equation}
\mathbf{P}^n= \mathbf{V}^\mathrm{T} \mathbf{p}^n
\,,\qquad 
   \mathbf{p}^n= \mathbf{V} \mathbf{P}^n  
\end{equation}
and we note that $\mathbf{V}\mathbf{V}^\mathrm{T}  = \mathbf{I}$, where $\mathbf{I}$ is an
identity matrix.  These matrix-vector products represent the Discrete Cosine Transform (DCT-II)
and its inverse~\cite{strang1999discrete}, respectively, which can be computed efficiently using
a FFT-based algorithm~\cite{vanloan1992computational}.  

The fully-discrete analogue to the modal update~\eqref{eqnmodalupdate} is now:
\begin{equation}
   \mathbf{P}^{n+1} = \mathbf{A}\circ\mathbf{P}^{n} - \mathbf{B}\circ\mathbf{P}^{n-1}
\end{equation}
where ``$\circ$'' denotes an element-wise product, and
\begin{subequations}
\begin{align}
\mathbf{A} &= 2\left[e^{-\sigma T}\cos(\omega_0 T), \dots,e^{-\sigma T}\cos(\omega_{Q-1}T)\right]^{\mathrm{T}}
\\
\mathbf{B} &= \left[e^{-2\sigma T}, \dots,e^{-2\sigma T}\right]^{\mathrm{T}}
\end{align}
\end{subequations}

Next we consider the boundary condition~\eqref{eqnbc0}, which we would like to
implement in a ``hard'' manner analogous to~\eqref{eqnbc1}.  In order to derive this, let us
introduce an intermediate grid function $\tilde{p}_l^{n}$, for which
$\tilde{p}_l^{n}=p_l^{n}$ only when $l>0$.
Let us also define:
\begin{subequations}
\begin{align}
\tilde{\mathbf{p}}^n &= \left[\tilde{p}^n_0,\tilde{p}^n_1,\dots,\tilde{p}^n_{N_x-1}\right]^{\mathrm{T}}
\\
\tilde{\mathbf{P}}^n &= \mathbf{V}^\mathrm{T} \tilde{\mathbf{p}}^n
\\
\boldsymbol{\delta} &= \left[1,0,\dots,0\right]^{\mathrm{T}}
\end{align}
\end{subequations}
i.e., $\tilde{\mathbf{p}}^n$ is the vector form of $\tilde{p}_l^{n}$ and $\tilde{\mathbf{P}}^n$ its modal form, and
$\boldsymbol{\delta}$ is a Kronecker delta in vector form.

We would like to express $p_l^{n+1}$ as $\tilde{p}_l^{n+1}$ plus a correction, such that~\eqref{eqnbc1} is
satisfied.  Clearly, this would be as simple as setting:
\begin{equation}
   p_0^{n+1} = \tilde{p}_0^{n+1} + (u^{n+1}-\tilde{p}_0^{n+1})
   \label{eqncorrectionFDTD}
\end{equation}
This could also be written as:
\begin{equation}
   \mathbf{p}^{n+1} = \tilde{\mathbf{p}}^{n+1} + \boldsymbol{\delta}(u^{n+1}-\boldsymbol{\delta}^{\mathrm{T}}\tilde{\mathbf{p}}^{n+1})
   \label{eqncorrect0}
\end{equation}
Left-multiplying by $\mathbf{V}^\mathrm{T}$, we have in the 
modal domain:
\begin{equation}
   \mathbf{P}^{n+1} = \tilde{\mathbf{P}}^{n+1} + \mathbf{V}^{\mathrm{T}}\boldsymbol{\delta}\left(u^{n+1} - \boldsymbol{\delta}^{\mathrm{T}}\mathbf{V}\tilde{\mathbf{P}}^{n+1}\right)
\end{equation}
Let us now define:
\begin{equation}
\boldsymbol{\phi} = \mathbf{V}^{\mathrm{T}}\boldsymbol{\delta} 
   = \left[\Phi_{0,0},\Phi_{1,0},\dots,\Phi_{Q-1,0} \right]^{\mathrm{T}}
\end{equation}
We can then propose the following scheme update with a ``hard source'' correction (forced-boundary):
\begin{subequations}
   \label{eqnmodalupdatematrix}
\begin{align}
   \tilde{\mathbf{P}}^{n+1} &= \mathbf{A}\circ\mathbf{P}^{n} - \mathbf{B}\circ\mathbf{P}^{n-1}
\\
   \mathbf{P}^{n+1} &= \tilde{\mathbf{P}}^{n+1} + \boldsymbol{\phi}\left(u^{n+1} - \boldsymbol{\phi}^{\mathrm{T}}\tilde{\mathbf{P}}^{n+1}\right)
\end{align}
\end{subequations}
This update is exact for the lossless case since it is, in fact, the exact FDTD scheme
(\eqref{eqnfdtd} and~\eqref{eqnbc1} via~\eqref{eqncorrectionFDTD})
transformed to a modal update.\footnote{This will be not be shown, but the proof is
carried out using the known eigenvalues/eigenvectors of the corresponding Laplacian
matrix (see, e.g.,~\cite{strang1999discrete}) or by Taylor series expansion (see,
e.g.,~\cite{etgen1989accurate}).}  It is not expected to be exact in the presence of loss since
the hard boundary is applied at $x=X/2$ and not $x=0$, which means there is a small distance
over which dissipation is not incurred (this error is negligible for practical purposes).

Alternatively, one can apply a simpler ``soft source'' update, for which the modal update becomes:
\begin{equation}
\label{eqnmodalsoft}
   \mathbf{P}^{n+1} = \mathbf{A}\circ\mathbf{P}^{n} - \mathbf{B}\circ\mathbf{P}^{n-1} 
+ F^{n}\boldsymbol{\phi}
\end{equation}
where $F^{n}$ is the forcing term~\eqref{eqnforce}.\footnote{It may be noted that~\eqref{eqnmodalsoft} is similar to the modal time-stepping update
presented in~\cite{botts2015extension} -- particular in the choices of time-stepping update and of a cosine
modal basis.  The numerical modal method in~\cite{botts2015extension}, paired with domain decomposition, was
ultimately used for 3-D wave simulations with damping that similarly mimics air attenuation.
However, the scheme
here is limited to one spatial dimension, and the 
``soft'' forcing term given here is of a different nature than that in~\cite{botts2015extension}
-- derived here specifically to satisfy the boundary condition~\eqref{eqnbc0}.}  
This update is also exact in the lossless
case, but it is expected to be less accurate than~\eqref{eqnmodalupdatematrix} when $\sigma>0$.
However, this update has the advantage that stability is straightforward to establish through
Z-transform analysis.

\begin{figure}[t]
   \begin{subfigure}[b]{0.48\textwidth}
      \centering
      \includegraphics[scale=0.53,clip,trim=1cm 0cm 1cm 0cm]{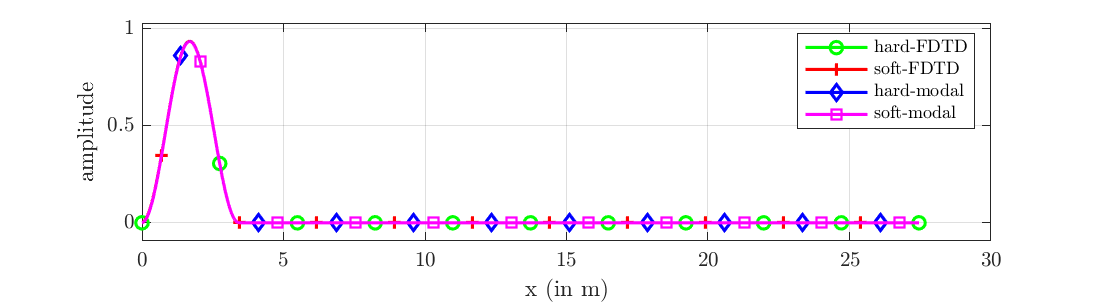}
      \caption{$t=10$\,ms}
   \end{subfigure}
   \begin{subfigure}[b]{0.48\textwidth}
      \centering
      \includegraphics[scale=0.53,clip,trim=1cm 0cm 1cm 0cm]{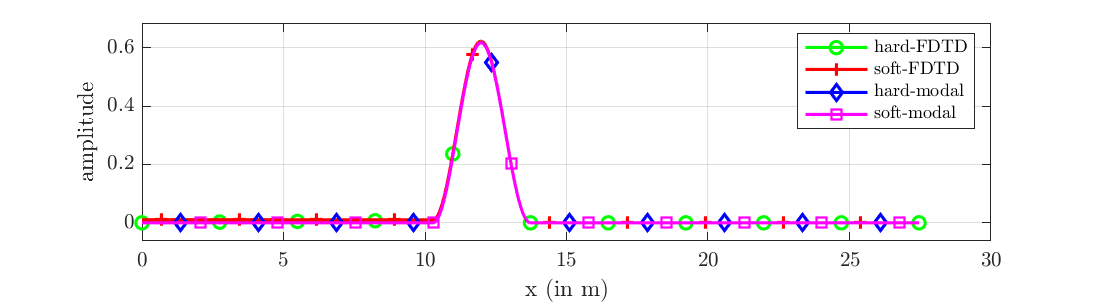}
      \caption{$t=40$\,ms}
   \end{subfigure}
   \begin{subfigure}[b]{0.48\textwidth}
      \centering
      \includegraphics[scale=0.53,clip,trim=1cm 0cm 1cm 0cm]{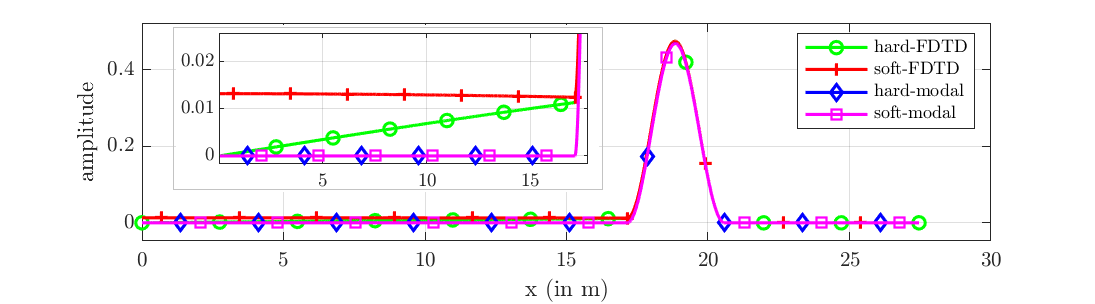}
      \caption{$t=60$\,ms}
   \label{figsim1c}
   \end{subfigure}
   \caption{Examples of simulated one-way dissipative wave propagation at different time
   instances $t$, using raised cosine input $u(t)$ with 60dB decay time
   $6\ln(10)\sigma^{-1}=1.0$\,s and speed of sound $c=343$\,m/s.  All schemes using 
   time-step $T=(48\textrm{\,kHz})^{-1}$.  The last image displays an overlay to show detail
   in wake behaviour.}
   \label{figsim1}
\end{figure}

\begin{figure*}[t]
   \begin{subfigure}[b]{0.24\textwidth}
      \centering
      \includegraphics[scale=0.40,clip,trim=0cm 0cm 1cm 0cm]{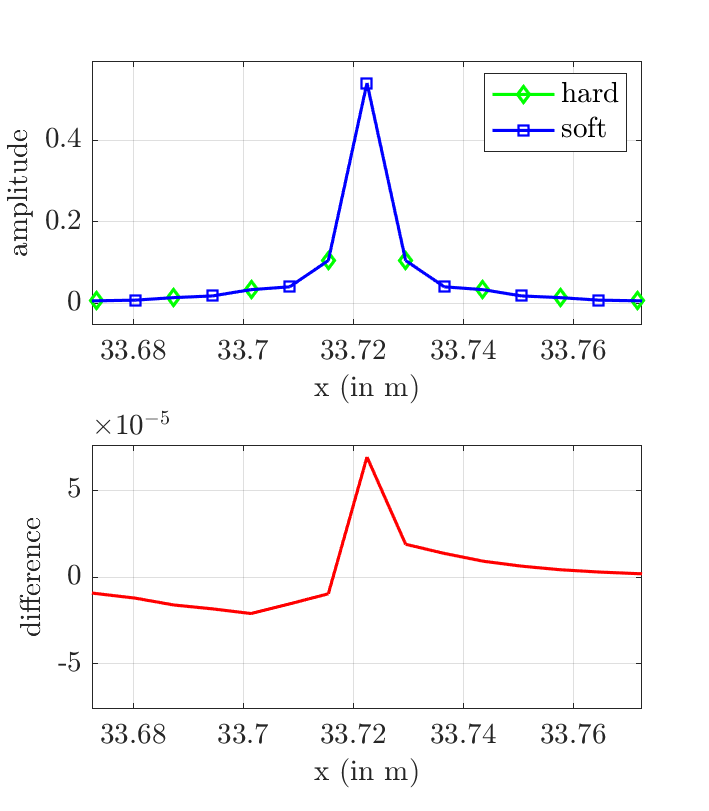}
      \caption{$t=0.1$\,s}
   \end{subfigure}
   \begin{subfigure}[b]{0.24\textwidth}
      \centering
      \includegraphics[scale=0.40,clip,trim=0cm 0cm 1cm 0cm]{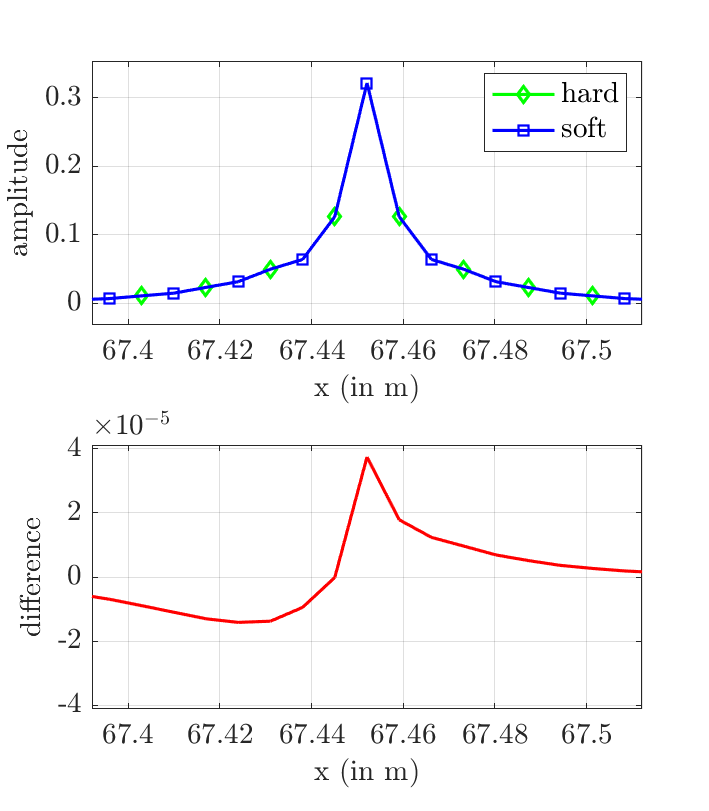}
      \caption{$t=0.2$\,s}
   \end{subfigure}
   \begin{subfigure}[b]{0.24\textwidth}
      \centering
      \includegraphics[scale=0.40,clip,trim=0cm 0cm 1cm 0cm]{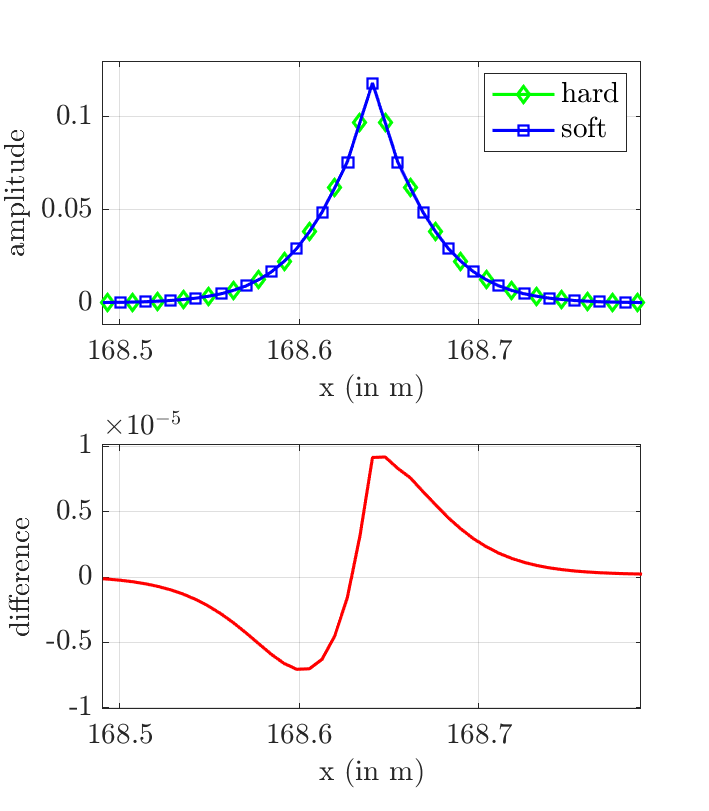}
      \caption{$t=0.5$\,s}
   \end{subfigure}
   \begin{subfigure}[b]{0.24\textwidth}
      \centering
      \includegraphics[scale=0.40,clip,trim=0cm 0cm 1cm 0cm]{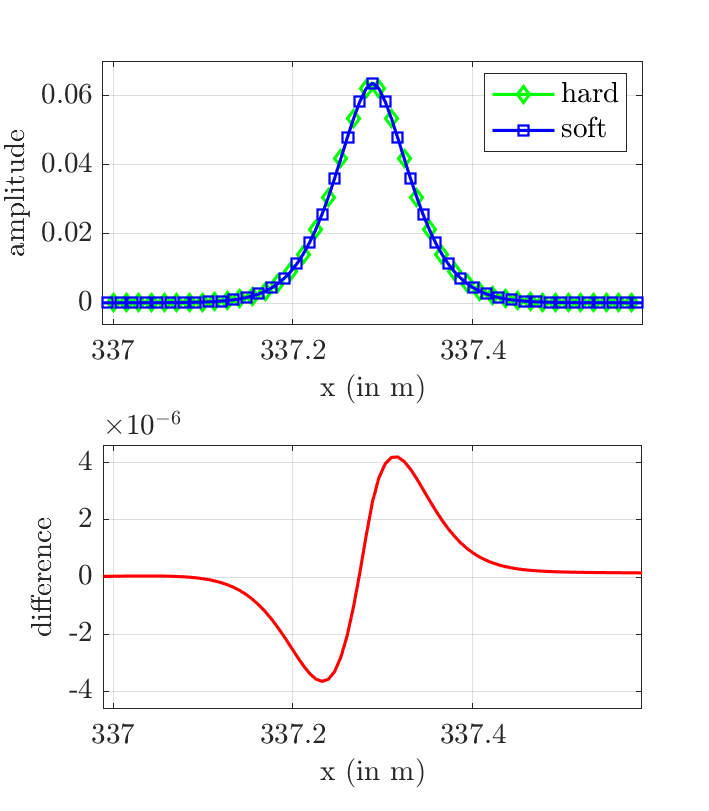}
      \caption{$t=1.0$\,s}
   \end{subfigure}
   \caption{Examples of simulated one-way dissipative wave propagation with modal schemes using
   ``soft'' and ``hard'' boundary conditions (and differences), at different time
   instances $t$, using Kronecker delta input ($u^n=[1,0,\dots,0]$) with decay constants
   $\sigma_q$ set from air
   attenuation model for air temperature $10^\circ$C and 20\% relative humidity ($c=337$\,m/s).
   Modal time-stepping using time-step $T=(48\textrm{\,kHz})^{-1}$.}
   \label{figsim2}
\end{figure*}
\subsection{Validation of frequency-independent dissipative schemes}
A simple test is conducted to evaluate the numerical schemes presented thus far.  The case of
zero loss will not be tested (for brevity) since exactness may be verified by hand.  As for
the lossy case $\sigma>0$,
Fig.~\ref{figsim1} displays simulated one-way dissipative wave propagation with a raised-cosine
input signal $u(t)$ with dissipation set for a 60dB decay time of $6\ln(10)\sigma^{-1}=1.0$\,s.
Displayed in the figure are $p(t,x)$ simulated using the FDTD schemes with ``hard'' and
``soft'' boundary updates, as well as using the corresponding modal schemes.  It can be seen
that the four simulations display similar behaviour, but differences are apparent in the wake of
the travelling wave (see zoom overlay in Fig.~\ref{figsim1c}).   For example, the ``soft-FDTD''
does not satisfy the boundary condition with the chosen raised-cosine input
(one can note a slight DC offset at $x=0$), whereas the ``hard-FDTD'' output displays a slope in
the wake of the wave which is absent from modal outputs.  In regards to the ``soft'' updates, it
should be noted that in this case the assumption $\sigma\ll \omega$ is not satisfied for low
frequencies.

For the remainder of this paper we will consider only modal methods, but it is worth mentioning
that what follows could be achieved with FDTD methods using a more general -- and
indeed more ``physically valid'' -- second-order wave equation with air absorption processes
included (after~\cite{hamilton2020viscothermal}).  However, in that case we would be faced with
numerical dissipation errors -- even with a domain limited to one spatial
dimension\footnote{Numerical phase velocity errors will be minimal in such 1-D FDTD schemes (and
non-existent in the lossless case~\cite{hildebrand1968finitedifference,ames1977numerical}), but
numerical \emph{dissipation} errors may be significant in high frequencies when discretising loss terms using
conventional means (trapezoid rule) -- e.g., as reported
in~\cite{hamilton2014visc,hamilton2020viscothermal} with regards to the problem of air
absorption.  Without resorting to higher-accurate schemes (with less favourable stability
properties), such dissipation errors can only be overcome with oversampling, which would
increase computational costs.  An analogous issue has been reported in the context of FDTD
modelling of musical strings~\cite{desvages2018physical}, where a possible solution is to
optimise parameters for fictitious loss terms in the discrete domain to minimize numerical
dissipation errors~\cite{desvages2019optimised}.} -- and this would necessitate grid refinement;
whereas the ``exact'' modal schemes are essentially free from dissipation errors (and are,
accordingly, more efficient for this specific application).

\section{Modal scheme for simulating air absorption\label{secmodal}}
In this section we extend the modal scheme to the frequency-dependent case to incorporate air
attenuation.  Since each mode in
the system is independent the modal damping coefficients can be set individually for each
mode (as in, e.g.,~\cite{botts2015extension,ducceschi2016plate}).  It suffices then to replace
$\sigma$ by $\sigma_q=c\alpha(\omega_q)$.  Note that with realistic air absorption (see, e.g.,
Fig.~\ref{figabs}), it can be assumed that $\alpha(\omega)\ll \omega/c$ for all frequencies of
interest (audible frequencies).

In terms of a PDE system, we have then $Q$ lossy wave equations ($q=0,\dots,Q-1$), each with its own damping
coefficient:
\begin{subequations}
   \label{eqnsyseries}
   \begin{align}
   \partial_t^2 p_q  + 2\sigma_q \partial_t p_q - c^2 \partial_x^2 p_q = 0
   \\
   p_q(t,0) = u_q(t)
   \end{align}
\end{subequations}
We take the associated
solution:
\begin{equation}
   p(t,x) = \sum_{q=0}^{Q-1} p_q(t,x) = \sum_{q=0}^{Q-1} \hat{u}_q e^{-\alpha_q x}e^{i(k_qx-\omega_q t)}
\end{equation}
where $ck_q= \sqrt{\omega_q^2 + \sigma_q^2}$, and 
where $u(t)$ is related to $\hat{u}_q e^{-i\omega_q t}$ through~\eqref{eqnufourier}.

Our fully-discrete modal scheme with ``hard'' boundary keeps the form~\eqref{eqnmodalupdatematrix},
but now with
\begin{subequations}
\begin{align}
\mathbf{A} &= \left[2e^{-\sigma_0 T}\cos(\omega_0 T), \dots,2e^{-\sigma_{Q-1} T}\cos(\omega_{Q-1}T)\right]^{\mathrm{T}}
\\
\mathbf{B} &= \left[e^{-2\sigma_0 T}, \dots,e^{-2\sigma_{Q-1} T}\right]^{\mathrm{T}}
\end{align}
\end{subequations}
Meanwhile, the modal scheme with the ``soft source'' boundary now takes the form:
\begin{equation}
\label{eqnmodalsoftq}
   \mathbf{P}^{n+1} = \mathbf{A}\circ\mathbf{P}^{n} - \mathbf{B}\circ\mathbf{P}^{n-1} 
+\mathbf{F}^{n}\circ\boldsymbol{\phi}
\end{equation}
with $\mathbf{A}$ and $\mathbf{B}$ changed as noted previously, and the forcing term taking the
\emph{vector} form:  
\begin{subequations}
\begin{align}
\mathbf{F}^n &= \left[F^n_0,F^n_1,\dots,F^n_{Q-1}\right]^{\mathrm{T}}
\\
F_q^n&=\frac{1+\sigma_q T/2}{1+\sigma_q T} u^{n+1} - \frac{1-\sigma_q T/2}{1+\sigma_q T} u^n
\end{align}
\end{subequations}

\subsection{Validation of modal schemes and comparisons with filter methods}
Numerical tests are presented in order to compare and validate the frequency-dependent modal
schemes in the context of air absorption.  For the following tests, air conditions of
$10^\circ$C temperature and 20\% relative humidity are chosen so that air absorption
deviates significantly from classical power-law attenuation in the audible frequencies (see
Fig.~\ref{figabs}).

Fig.~\ref{figsim2} shows simulated one-way dissipative wave propagation with a
Kronecker delta input with dissipation set to mimic the desired air attenuation. Displayed in
the figure are $p(t,x)$ simulated using the modal schemes with ``hard'' and ``soft'' boundary
update, and differences between them.  It can be seen that the schemes return similar results
with small differences.  



\begin{figure*}[ht]
   \centering
   \begin{subfigure}[b]{0.32\textwidth}
      \includegraphics[scale=0.34,clip,trim=0cm 0cm 1cm 0cm]{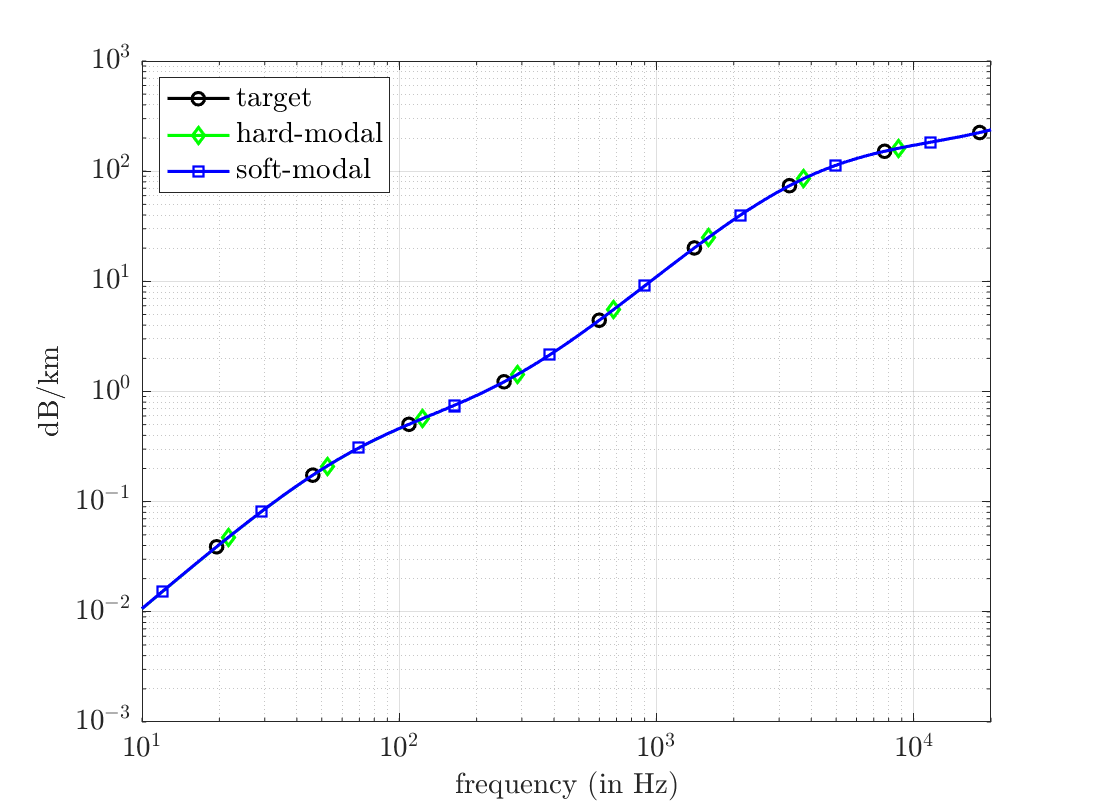}
      \caption{Modal approach}
      \label{figsim3a}
   \end{subfigure}
   \begin{subfigure}[b]{0.32\textwidth}
      \includegraphics[scale=0.34,clip,trim=0cm 0cm 1cm 0cm]{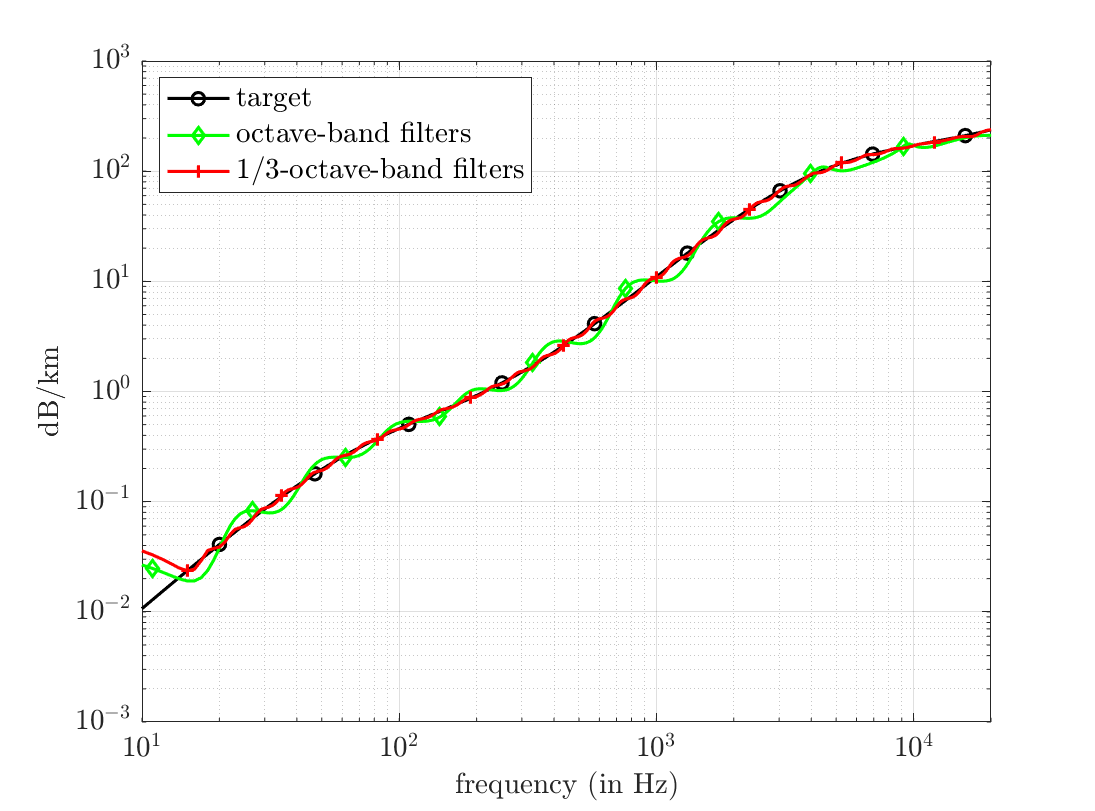}
      \caption{Octave-band filter approach}
      \label{figsim3b}
   \end{subfigure}
   \begin{subfigure}[b]{0.32\textwidth}
      \includegraphics[scale=0.34,clip,trim=0cm 0cm 1cm 0cm]{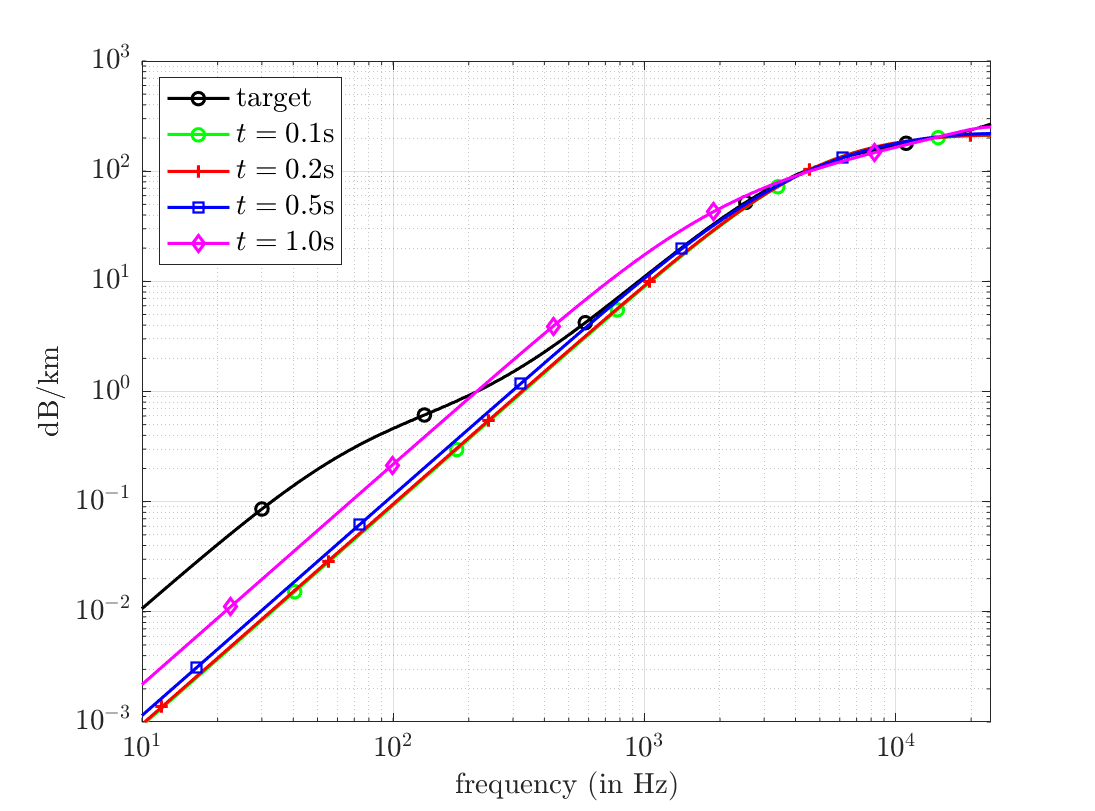}
      \caption{Time-varying optimised lowpass-filter approach}
      \label{figsim3c}
   \end{subfigure}
   \caption{Simulated air attenuation using modal schemes and digital-filter approaches, and
   target air attenuation for air temperature $10^\circ$C and 20\% relative humidity
   ($c=337$\,m/s).  Modal time-stepping using time-step $T=(48\textrm{\,kHz})^{-1}$.  Filter
   methods use 48\,kHz sample rates.  Octave-band filters are of third-order Butterworth type,
   and optimised lowpass-filter approach is using the method presented in~\cite{kates2020airabsorption}.}
\end{figure*}

To validate the modal schemes in their reproduction of a target air attenuation, we can simply
measure the attenuation at a given time instant (corresponding to a distance travelled) and
compare to the ``analytic'' air absorption model (from~\cite{iso19939613}).  For the same
conditions above, Fig.~\ref{figsim3a} shows the attenuation per distance (in km) simulated,
calculated at time instant $t=1.0$\,s, along with the target air attenuation curve.  It can be
seen that agreement with the target air attenuation is excellent (as expected for these modal
schemes).   Results are identical for other time instances and are thus left out for
brevity.\footnote{However, for very large distances dissipation can be sufficient to exceed the
relative accuracy of finite-precision, and -- only in those cases -- simulated measurements of
air attenuation will not return accurate results.}

In contrast, a more conventional filter approach is plotted in Fig.~\ref{figsim3b}, in which
third-order Butterworth octave-band filters are used as a analysis/synthesis filter bank used to
process a 48\,kHz-sample-rate Kronecker delta input with distance-based exponential damping applied to mimic air
absorption (sampled at 11 octave-band centre frequencies) -- as may be found in traditional
geometrical acoustics room simulations such as, e.g.,~\cite{scheibler2018pyroomacoustics}.
Ripples are apparent in the attenuation simulated by the filterbank approach, as one might
expect.  This can be improved using a filterbank of 1/3-octave band filters (as in,
e.g.,~\cite{schroder2011physically}), which, if chosen well, leads to more tightly-fitting
ripples, as seen in the same figure (here, using 32 1/3-octave-band filters of third-order
Butterworth type).  Another potential option would be the use of FIR partition-of-unity
filters~\cite{antoni2010orthogonal}.

Additionally, results from an alternative, recently-presented~\cite{kates2020airabsorption},
optimised IIR filter approach is plotted in Fig.~\ref{figsim3c}.  This method is based on a set of triple-cascaded
first-order lowpass filters optimised for attenuation tuned to distances corresponding to each sampling
instant~\cite{kates2020airabsorption}.  A Kronecker delta is filtered to measure
simulated attenuation tuned to different time instants (distances).  In the figure it can be seen that this approach does well to reproduce
attenuation in high frequencies (where it is most significant), but errors appear in low
frequencies.  On the other hand, one can note that simulated attenuation curves vary as a
function of time instant, and for longer times (higher attenuation) attenuation
filters deviate more from the target attenuation (in dB/km).\footnote{To produce data for the
optimised lowpass filter approach~~\cite{kates2020airabsorption}, Matlab code provided as
supplementary material to~\cite{kates2020airabsorption} was used, but reconfigured to optimise over
20\,Hz--20\,kHz with approximately 40k frequency samples.}  

\subsection{Adding air absorption to simulated RIRs}
It remains to detail the procedure to add air absorption to a
pre-simulated RIR (simulated without air absorption) using the presented modal schemes.  Consider then a finite-duration discrete
RIR signal $h[n]$ with sample rate $1/T$ (e.g., 48\,kHz), with $N_t$ samples, indexed by
$n'=0,\dots,N_t-1$.   Of course, $h[n]$ can be seen as a linear combination of $N_t$ shifted and
scaled Kronecker deltas.  We have shown that the distance-based attenuation is accurately
reproduced for a Kronecker delta input, so it follows that we can simply apply the same procedure to
the signal itself to process each individual sample.  Since the system is linear and
time-invariant, all such samples can be processed simultaneously (i.e. in a single modal
simulation).  

We must ensure that correct distance attenuation is applied to each sample
$n$, where that distance would be simply $d=c(n+1/2)T$ (in m).\footnote{Here we assume a constant
speed of sound $c$, but it is important in practice that any ``pre-delay'' in the signal be
removed (or compensated for in the application of attenuation).}  For this we must \emph{time-reverse} the
input RIR as we feed it into the system (via the left-boundary condition);
i.e., we can set the forcing signal to our ODE system to be:
\begin{equation}
u^n = h[N_t-1-n]
\end{equation}
We can then run the forced modal schemes with $Q=N_t$ modes for $N_t$ time-steps, starting from zero initial conditions
($\mathbf{P}^{-1}=\mathbf{P}^{-2}=0$).  The output RIR with air attenuation added, which we
will denote $h'[n]$, can then be read from the final state of our wave simulation as follows:  
\begin{equation}
h'[n] = p^{N_t-1}_n\,, \quad n=0,\dots,N_t-1
\end{equation}
where $p^{N_t-1}_n$ are simply the elements of $\mathbf{p}^{N_t-1} = \mathbf{V}\mathbf{P}^{N_t-1}$.
 Note that
the entire simulation can be run in the modal domain.  Projection back to a spatial grid (with
$N_x$ points) is only needed at the end, and this can be accomplished using a fast inverse DCT-II
algorithm~\cite{vanloan1992computational}.

When losses are disabled ($\sigma_q=0$), the proposed procedure will simulate the solution~\eqref{eqnbc1} and
return $h'[n]=h[n]$, as desired.  Also, it is recommended that the ``soft'' update be used for the
application of air absorption, since we have seen that ``hard'' and ``soft'' boundary
updates give similar results, and the soft update ultimately requires fewer operations.

\subsection{A reference solution}
At this point it is useful to introduce a ``reference solution'' for the problem of adding air
absorption to a RIR.  Taking $h[n]$ and $h'[n]$ as a $N_t\times 1$ vectors:
\begin{subequations}
   \begin{align}
      \mathbf{h} &= \left[h[0],h[1],\dots,h[N_t-1]\right]^{\mathrm{T}}
      \\
      \mathbf{h}' &= \left[h'[0],h'[1],\dots,h'[N_t-1]\right]^{\mathrm{T}}
   \end{align}
\end{subequations}
We can propose a reference solution:
\begin{equation}
   \mathbf{h}' = \boldsymbol{\Omega}\mathbf{V}^\mathrm{T} \mathbf{h}
   \label{eqnrefsol}
\end{equation}
where $\boldsymbol{\Omega}$ is a $N_t \times N_t$ matrix with elements at rows $q$ and columns
$l$ (with zero-based indexing) given by:
\begin{equation}
   \Omega_{q,l} = e^{-\alpha_q x_l} \Phi_{q,l}
\end{equation}
where $x_l = (l+1/2)X$ and where $\Phi_{q,l}$ is given by~\eqref{eqncosinePhi}.  This reference
solution is simply a decomposition of a signal onto an orthonormal cosine basis and a subsequent
projection back onto the same basis weighted by distance-based air attenuation, under the
assumption that the $l$th sample travels $x_l$ metres.  This can ultimately be seen as a
limiting case of the frequency-band filterbank approach, where here $N_t$ frequency bands are
used.  Note that when $\alpha_q=0$~\eqref{eqnrefsol}, $\boldsymbol{\Omega}=\mathbf{V}$,
and~\eqref{eqnrefsol} returns $\mathbf{h}'=\mathbf{h}$ (as desired).   However, this is not
claimed to be an analytic solution to the dissipative wave problems seen before, nor to the more
general second-order wave systems with air absorption
processes~\cite{pierce2019acousticsChap9,hamilton2020viscothermal}; it should be seen as a
``brute-force'' application of air attenuation under the assumptions used in this study.  

In terms of its compute cost, for $\alpha_q>0$ this reference solution is somewhat
impractical for the reason that $N_t$ should on the order of $10^4$--$10^5$, and only the
matrix-vector product $\mathbf{y}=\mathbf{V}^\mathrm{T}\mathbf{h}$ (DCT-II) can use a fast
FFT-based algorithm~\cite{vanloan1992computational}.  The matrix-vector product
$\boldsymbol{\Omega}\mathbf{y}$ must be computed directly, and this would require $N_t^2$
evaluations of transcendental functions, and $N_t^2$ storage (if $\boldsymbol{\Omega}$ is
precomputed).  More will be said about computational complexity later.

\subsection{Numerical examples}
\begin{figure}[t]
   \begin{subfigure}[b]{0.48\textwidth}
      \centering
      \includegraphics[scale=0.48,clip,trim=0.4cm 0cm 0.5cm 0cm]{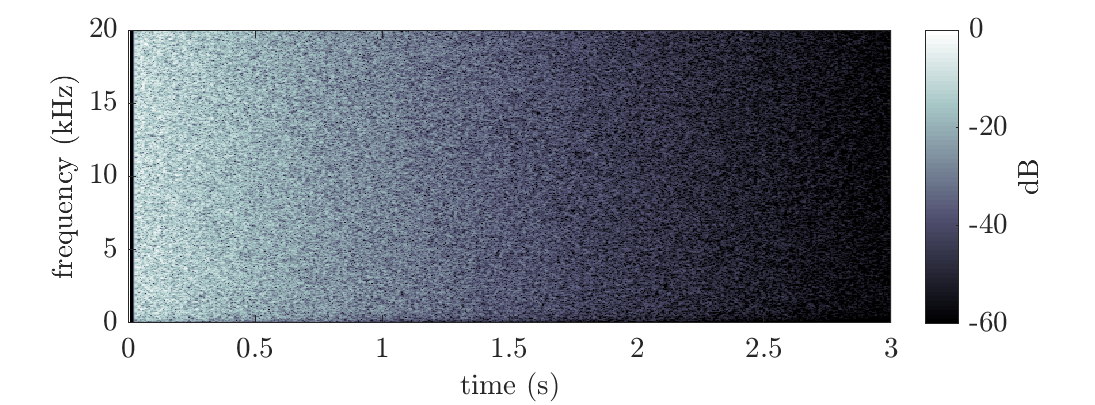}
      \caption{Spectrogram of $h[n]$, obtained with lossless-air image-source method}
      \label{figimage1}
   \end{subfigure}
   \begin{subfigure}[b]{0.48\textwidth}
      \centering
      \includegraphics[scale=0.48,clip,trim=0.4cm 0cm 0.5cm 0cm]{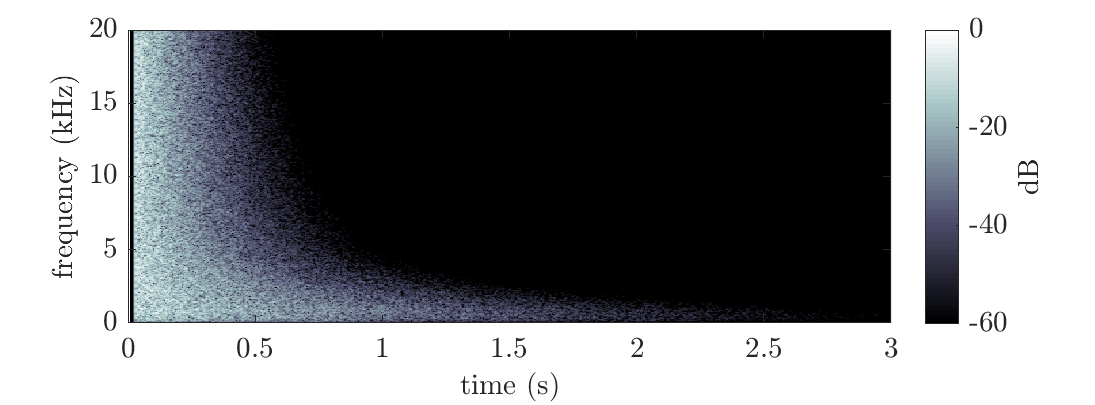}
      \caption{Air attenuation added via modal scheme}
      \label{figimage2}
   \end{subfigure}
   \begin{subfigure}[b]{0.48\textwidth}
      \centering
      \includegraphics[scale=0.48,clip,trim=0.4cm 0cm 0.5cm 0cm]{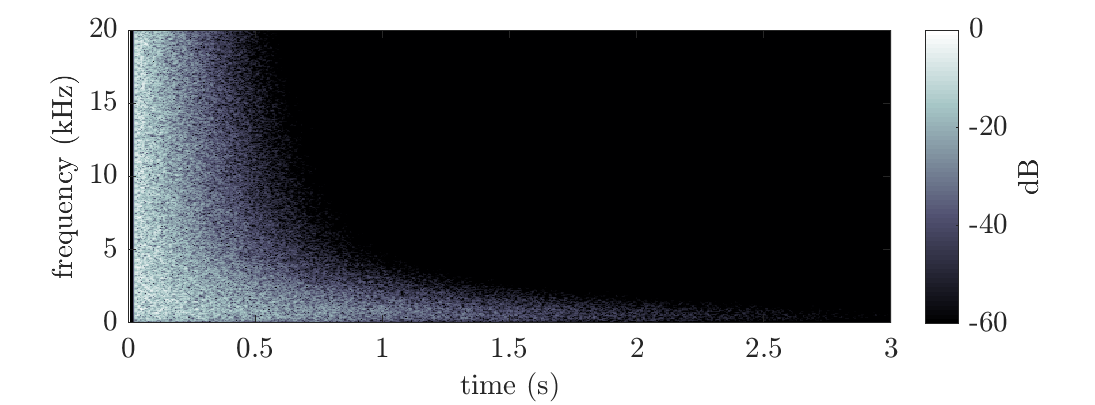}
      \caption{Air attenuation added with reference solution}
      \label{figimage3}
   \end{subfigure}
   \begin{subfigure}[b]{0.48\textwidth}
      \centering
      \includegraphics[scale=0.48,clip,trim=0.4cm 0cm 0.5cm 0cm]{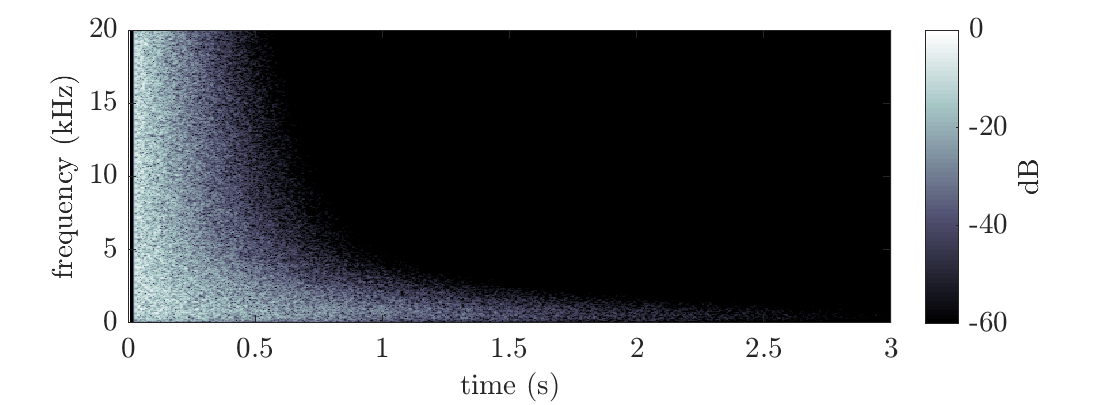}
      \caption{Air attenuation added via STFT processing}
      \label{figimage4}
   \end{subfigure}
   \caption{Spectrograms of room impulse response, simulated with image source method, with and
   without air absorption added.}
   \label{figISM}
\end{figure}

In order to test the application of air absorption to a RIR using this modal scheme (with
``soft'' boundary update), we conduct a simple image source simulation~\cite{allen1979image}.  A
room of size $7.2\times5.1\times4.3$ (in m) is chosen with a Sabine-absorption coefficient of
$0.045$, with a source and receiver at positions $(6.1,2.0,1.1)$ and $(3.3,3.1,1.3)$
respectively.  Image source distances are randomly perturbed by a relative amount in
range $[-1\%,+1\%]$ in order to eliminate sweeping echoes due to perfect rectangular
symmetry~\cite{desena2015modeling} (and ultimately increasing diffuseness in the room).  Air
conditions are as before (10$^\circ$C and 20\% relative humidity) with a corresponding speed of
sound $c=337$\,m/s.   A RIR, $h[n]$, is simulated at 48kHz, and its spectrogram -- computed
with 1024-sample Hann-windowing and 75\% frame overlap -- is displayed in Fig.~\ref{figimage1}.
Using $h[n]$ as input to the modal scheme (using the time-reversed procedure described
previously) returns a modified RIR, $h'[n]$, whose spectrogram is displayed in
Fig.~\ref{figimage2}, where the effect of air absorption is clearly seen.  The corresponding reference
solution is shown in Fig.~\ref{figimage3}.  Additionally, air attenuation is added to $h[n]$ with the
overlap-add (OLA) approach (STFT processing) (after~\cite{southern2013room,saarelma2018challenges})
and shown in Fig.~\ref{figimage4}.\footnote{For OLA processing, a Hann window of length 1024
samples was used for analysis and synthesis windows, scaled for a constant unity overlap-add
envelope (with sufficient zero padding to provide perfect reconstruction in lossless case), with
air absorption applied to STFT frames by frequency-domain multiplication.}  Agreements are
excellent between the three methods, at least from the point of view of spectrograms. 

\begin{figure}[t]
   \centering
   \includegraphics[scale=0.50,clip,trim=0.5cm 0cm 1cm 0cm]{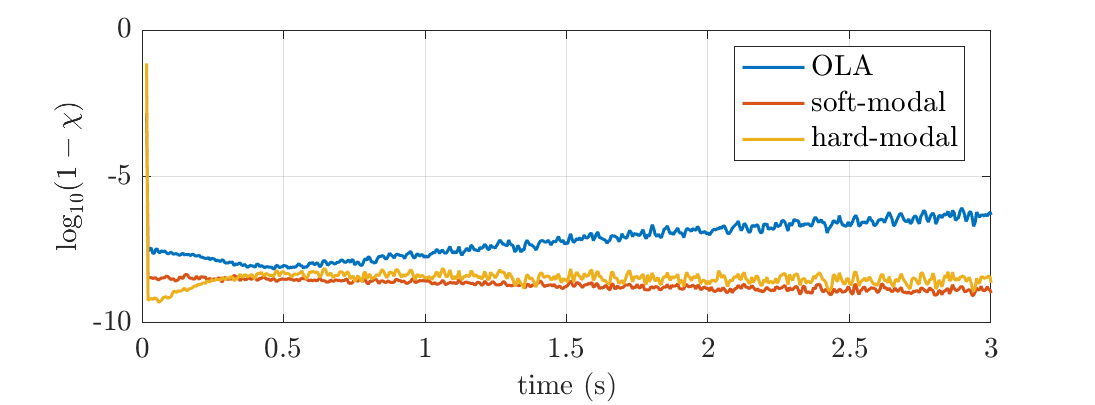}
   \caption{RIR STFT frame similarities~\eqref{eqnsimilarity} to reference solution (see Fig.~\ref{figISM}).}
   \label{figismsim}
\end{figure}

For a more objective comparison, it is useful to compute similarities between STFT frames
corresponding to the spectrograms in Fig.~\ref{figISM}.  This can be accomplished with the
following geometric similarity index: 
\begin{equation}
   \chi[n] = \frac{\| Y_1(n,\omega)Y_2^*(n,\omega)\|_{\omega}}{\| Y_1(n,\omega)\|_{\omega} \| Y_2(n,\omega)\|_{\omega}}
   \label{eqnsimilarity}
\end{equation}
Here $Y_1(n,\omega)$ and $Y_2(n,\omega)$ represent the STFTs of two discrete-time signals $y_1[n]$ and
$y_2[n]$, and the Euclidean norms are calculated across
the frequency-dimension, and $0\leq \chi[n]\leq 1$.  Comparing to the reference solution, in Fig.~\ref{figismsim} we plot
$\log_{10}(1-\chi)$ for the hard-modal output, and the output from the OLA method, along with the output of the soft-modal
scheme (corresponding spectrogram not shown for brevity).  It can be seen that the outputs generated are extremely
similar to the reference solution, with modal schemes performing best, generally by an order of
magnitude or more (in terms of frame similarity).  Indeed, the OLA
method, at least as implemented here, would seem to perform more than well enough for practical use.
With that said, it should be remarked that this constitutes the first evaluation of the
accuracy of the OLA approach in the literature (to the author's knowledge).

\section{Conclusions and Final Remarks\label{secconc}}
A method to add general air absorption to simulated room impulse responses was provided in this
paper.
The method uses a fully-discrete modal scheme with boundary conditions chosen for one-way 
dissipative wave propagation, such that accurate frequency-dependent, distance-based air
attenuation may be applied to an input RIR.  FDTD schemes were used to help derive a ``soft
source'' boundary update, which was also compared to a ``hard'' forced-boundary update.
Frequency-dependent modal schemes were validated against target air attenuation curves, along
with existing filter-based approaches.  An image source room simulation was used to test the
application of air attenuation to a RIR, and this was compared to air attenuation applied via
STFT-based time-frequency processing.


The method presented is advantageous for its generality (it can apply air attenuation for any air
conditions) and because it does not need tuning or optimisation as required in filter-based methods.
On the other hand, this method is more computationally demanding than filter-based methods,
which would be more applicable to real-time applications.  This modal approach is better suited to
offline data generation (as needed for, e.g.,  auditory research~\cite{fogerty2020effect} or
speech dereverberation~\cite{naylor2010speech}).  Additionally, this approach is well-suited to complement
(offline) RIR simulation with 3-D wave-based methods~\cite{hamilton2016gpuISMRA,lai2020ctk},
where modelling air absorption processes directly leads to significant increases in memory usage
(e.g., at minimum 50\% more memory than lossless wave-equation
schemes~\cite{hamilton2014visc,hamilton2020viscothermal}).

With regards to computational costs, the complexity of the presented modal is
$O(N_t^2)$ where $N_t$ is the number of samples to be computed, which is similar to a 1-D FDTD
scheme (in the lossy case~\cite{bilbao2009numerical}).  On the other hand, time-frequency
processing via STFT would be $O(\frac{N_t}{H}M\log_2 M)$ where $M$ is the frame size (assumed
power of two) and $1\leq H\leq M$ is the hop size in samples.  While it is possible for the OLA
method to have higher complexity than the modal approach (depending on choices of $M$ and $H$),
in practice it is expected that the complexity should be lower for STFT-based processing (usually $M\ll
N_t$) .  Low-order IIR filter-based approaches would usually be $O(N_t)$, but frequency-band
approaches would reach $O(N_t^2)$ in the high-resolution limit of $N_t$ frequency bands.  The
reference solution provided here can be seen as such a limiting filterbank -- and is accordingly also
$O(N_t^2)$ -- but should be less efficient than modal recursion because of the need for evaluation of
$N_t^2$ exponentially-damped-cosine matrix elements. 

In regards to compute times, it must also be noted that modal updates are highly, and easily,
parallelisable (each mode can be updated in parallel), and processing a RIR with this modal
approach could be accomplished in only a few seconds on a modern desktop computer (with an
optimised multithreaded code -- see,
e.g.,~\cite{ducceschi2016plate}), even with $N_t$ on the order of $10^5$.  To keep
the number of modes ($Q=N_t$) to a minimum the sampling rate of the input RIR should be not be
higher than necessary for audio purposes.

Finally, it should be mentioned that the presented approach is only applicable to static
scenes (precluding moving sources or receivers), but allows for source/receiver
directivity~\cite{bilbao2018directional,bilbao2019local}.  It could also be used for outdoor
scenes provided that wind is not significant.  On the other hand, for typical indoor conditions
(e.g., 20$^\circ$C, 50\% humidity), it may be sufficient to view air absorption as a simple
power-law attenuation, in which case the Green's function method recently presented
in~\cite{hamilton2021dafx} may be preferable with its $O(N_t)$ complexity and applicability
to real-time usage.  In future work, the various methods of applying air attenuation to a
simulated RIR could be compared on a perceptual basis (through listening tests).



\section*{Acknowledgment}
Thanks to Michele Ducceschi, Charlotte Desvages, and James-Michael Leahy for fruitful
discussions on topics relating to numerical modelling of dissipative wave propagation.

\appendix 
\renewcommand{\theequation}{\thesection.\arabic{equation}}
\setcounter{equation}{0}
\section{Appendix}
\subsection{Boundary condition for lossy wave equation}
We consider a semi-infinite domain terminated at $x=0$ (and extending towards $+\infty$) with
PDE~\eqref{eqnlossywave} and desired boundary
condition~\eqref{eqnbc0}.

An energy analysis reveals that one can set the left boundary as:
\begin{equation}
   \pd_x p(t,x)\rvert_{x=0} = -f(t)
   \label{appeqnbc1}
\end{equation}
where $f(t)$ is an unspecified forcing signal to be associated to $u(t)$.
Following~\cite{kinsler2000fundamentals}, but here with dissipation included, one can impose a
rightward-travelling plane-wave solution of the form:
\begin{equation}
   p(t,x) =  \hat{p}e^{i(\hat{k} x -\omega t)} 
\end{equation}
where $\hat{k}=k+i\alpha$.

By linearity, we take $u(t)=\hat{u}e^{-i\omega t}$ and $f(t)=\hat{f}e^{-i\omega
t}$.  From~\eqref{eqnbc0} and~\eqref{appeqnbc1}, we get at the left boundary:
\begin{equation}
   p(t,0) = \hat{u}e^{-i\omega t}
\end{equation}
and
\begin{equation}
   \pd_x p(t,x)\rvert_{x=0} = i(k+i\alpha) u(t)
\end{equation}
which implies $f(t) = -i\hat{k} u(t)$.
With $k\approx \omega /c$, we then have:
\begin{equation}
   \pd_x p(t,x)\rvert_{x=0} = -\frac{1}{c}(c\alpha - i\omega)u(t)
\end{equation}
and
\begin{equation}
   \label{appeqnbc3}
   \pd_x p(t,x)\rvert_{x=0} = -\frac{1}{c}(\sigma u(t) + \dot{u}(t))
\end{equation}

\subsection{Soft-source FDTD boundary condition for lossy wave equation}
We can discretise~\eqref{appeqnbc3} at $l=0$ with:
\begin{equation}
   \frac{p^n_{0}-p^n_{-1}}{X} = -\frac{1}{c}\left(\sigma \frac{u^{n+1} + u^{n}}{2} + \frac{u^{n+1} - u^{n}}{T}\right)
\end{equation}
Solving for the ghost point $p^n_{-1}$ with $T=X/c$, and inserting into~\eqref{eqnfdtd} at $l=0$, we get:
\begin{multline}
   p^{n+1}_0 =  \frac{1}{1+\sigma T}\left(p^n_{1} + p^n_{0} +(\sigma T-1) p^{n-1}_{0}\right) \\
+ \underbrace{\frac{1+\sigma T/2}{1+\sigma T} u^{n+1} - \frac{1-\sigma T/2}{1+\sigma T} u^n}_{F^n}
\end{multline}
where the ``soft source'' forcing term is denoted $F^n$.
It is straightforward to check by hand that~\eqref{eqnwavedisc} satisfies this update when
$\sigma=0$, and thus remains exact in the lossless case.  Furthermore, in the lossless case, this can
be seen as a consistent discretisation of:
\begin{subequations}
   \begin{align}
      \pd_t^2 p - c^2 \pd_x^2 p = c\delta(x-\varepsilon)\dot{u}(t) &
      \\
      \pd_x p\rvert_{x=0} = 0&
   \end{align}
\end{subequations}
in the limit of $\varepsilon\to 0$ ($\varepsilon\geq 0$).  This is a system that has the identical one-way travelling wave
solution~\eqref{eqntrav}, which can be verified through the known Green's function to
the lossless 1-D wave equation and the principle of mirror images (with two travelling waves
coalescing in the limit of $\varepsilon\to 0$).  By the same token, the lossy FDTD scheme can be seen as
a discretisation of the following PDE system:
\begin{subequations}
   \label{eqnappforce}
   \begin{align}
      \pd_t^2 p + 2\sigma \pd_t p - c^2 \pd_x^2 p = c\delta(x-\varepsilon)(\sigma u(t) + \dot{u}(t))&
      \\
       \pd_x p\rvert_{x=0} = 0 &
      \label{eqnappforceb}
   \end{align}
\end{subequations}
in the limiting case of $\varepsilon\to 0$.  This would imply that the system consisting
of~\eqref{eqnlossywave} and~\eqref{appeqnbc3} has the same
solutions as~\eqref{eqnappforce}, under the assumption that $\frac{1}{X}\delta[0]$, a scaled
Kronecker delta, is a discretisation of $\delta(x-\varepsilon)$ in the limit of $\varepsilon\to 0$ on
the chosen grid (where $l=0$ pertains to $x=X/2$).  As such, the forcing term $F^n$ may also be used
in the modal update~\eqref{eqnmodalsoft}, where the chosen cosine basis already satisfies the boundary
condition~\eqref{eqnappforceb}.

\bibliographystyle{IEEEbib}
\bibliography{refs}

\end{document}